\documentclass[a4paper,11pt]{article}
\usepackage[margin=3cm]{geometry}
\usepackage{amsfonts,amsmath,amssymb,amscd,graphicx}
\usepackage[colorlinks]{hyperref}
\usepackage{mathrsfs}

\newtheorem{theorem}{Theorem}[section]
\newtheorem{corollary}[theorem]{Corollary}
\newtheorem{lemma}[theorem]{Lemma}

\newtheorem{proposition}[theorem]{Proposition}

\newtheorem{example}[theorem]{Example}
\numberwithin{equation}{section}

\newenvironment{preuve}[1][]
{\vskip 2mm  \noindent\emph{\bf Proof#1. }}{$\Box$ \vskip 2mm}

\newcommand{\beq}{\begin{eqnarray}}
\newcommand{\eeq}{\end{eqnarray}}
\newcommand{\beqe}{\begin{eqnarray*}}
	\newcommand{\eeqe}{\end{eqnarray*}}

\DeclareMathOperator{\vol}{vol}
\DeclareMathOperator{\inte}{Int}

\DeclareMathOperator{\area}{area}

\newcommand{\cov}{\ensuremath \mathrm{Cov}}

\newcommand{\End}{\ensuremath \mathrm{End}}
\newcommand{\id}{\ensuremath \mathrm{Id}}

\newcommand{\Nn}{\mathbb{N}}
\newcommand{\R}{\mathbb{R}}

\newcommand{\C}{\mathbb{C}}

\let\epsilon=\varepsilon

\begin{document}

\title{ How curved is a random complex curve?}
\author{\sc Michele Ancona and Damien Gayet}
\maketitle
\begin{abstract}
	In this paper, we study the curvature properties of random complex plane curves. We bound from below the probability that a uniform proportion of the area of a random complex degree $d$ plane curve has a curvature smaller than $-d/8$. Our lower bound is uniform, in the sense that it does not depend on $d$. We also provide uniform upper bounds for similar probabilities. 
 These results extend to random complex curves of projective surfaces equipped with an ample line bundle. This paper can be viewed as a sequel of~\cite{AG}, where other metric statistics were given. On a larger time scale, it joins the general program initiated in~\cite{SZ99} of understanding random complex hypersurfaces of projective manifolds.  
\end{abstract}

\begin{flushright}
{\it	This paper is dedicated to the memory of Steve Zelditch}
	\end{flushright}

\tableofcontents
\section{Introduction}

\subsection{The standard projective setting}

Smooth complex projective curves of the complex projective plane $\C P^2$ have the remarkable property that their topology depends only on the degree $d$ of a defining homogeneous polynomial. More precisely,  they all are compact connected Riemann surfaces of genus  $\frac{1}2(d-1)(d-2)$.
 When equipped with the
restriction of the Fubini-Study metric $g_{\mathrm{FS}}$, these complex curves become Riemannian real surfaces. 
By the Wirtinger theorem, their area depends only on their degree:
$$ \forall d\geq 1, \forall P\in \C^{hom}_d[X_0, X_1,X_2], \ \area_{g_{\mathrm{FS}|Z(P)}}(Z(P))= d,$$
where $Z(P):=\{P=0\}\subset \C P^2.$
However, all the other Riemannian quantities strongly depend on the complex curve.
This paper deals with the statistics of \emph{curvature} properties of such surfaces, when they are chosen at random, for a fixed large degree. 

Let us recall some deterministic facts about the curvature of plane complex curves. 
For any $P\in \C^{hom}_d[X_0,X_1, X_2]$,  and any $x\in Z(P)$, denote by $K(x)$ the Gaussian curvature for the restriction of $g_{\mathrm{FS}}$ on $Z(P)$. 	The curvature satisfies the upper bound $$K\leq 2\pi,$$ 
which is achieved by complex lines in $\C P^2$, see for instance~\cite{ness1977curvature}. 
The points where $K=2\pi$ are called \emph{inflexion} points. 
In fact, 
for any $d\geq 2$ and any generic polynomial $P\in \C^{hom}_d[X_0,X_1, X_2]$, 
$$ 
\# \{x\in Z(P), K(x)=2\pi \}= 3d(d-2).$$
In particular, the degree 2 curves are the only one for which the curvature is strictly less than $2\pi$. 
Moreover, 
$$\forall d\geq 2, \inf_{P\in \C^{hom}_d[X_0,X_1, X_2]} 
\inf_{x\in Z(P)} K(x)= -\infty.$$
Indeed, the smoothing the union of a degree $d-1$ smooth complex curve and a complex line 
provides a degree $d$ curve whose curvature at the $d-1$ intersection points becomes infinitely negative 
when the smoothing becomes smaller and smaller.

Finally, by the Gauss--Bonnet theorem and the genus formula, 
\begin{equation}\label{GB}
\int_{Z(P)} K (x)\mathrm{d}\vol(x) =2-(d-1)(d-2),
\end{equation}
so that the average of $K$ over the complex curve is asymptotic to $-d$.

In this paper, we study the curvature of these complex curves when the defining polynomial is taken at random. 
Let
\begin{equation}\label{mesure2} P= \sum_{i+j+k =d} a_{i,j,k}\sqrt{ \frac{(d+2)!}{2 i!j!k!}}
X_0^{i}X_1^jX_2^{k}, \end{equation}
be a degree $d$ random polynomial,
where $(a_I)_{|I|=d}\in \C P^{N_d-1}$
are random coefficients chosen uniformly on the projective space $\C P^{N_d-1}$ equipped with its Fubini-Study metric (the quotient metric induced by the standard metric on $\C^{N_d}$) and $N_d= \dim \C^{hom}_d[X_0,X_1, X_2]$. We can also choose $(a_I)_I$ in the standard sphere $\mathbb S^{2N_d-1}\subset \C^{N_d}$, or being independent complex standard Gaussian random variables. 
We denote by $\mu_d$ the measure associated to~(\ref{mesure2}). This measure is naturally associated to the unique $U(3)$-invariant Hermitian product on $\C^{hom}_d[X_0,X_1, X_2]$, see Example~\ref{bobmarley}.

In general, our results will estimate the proportion of a complex curve where its curvature is controlled. Hence, for any real surface $Z$ equipped with a metric $g$ with finite area and any subset $A\subset \R$, let us define 
\begin{equation}\label{kappa0}
\kappa(Z,g,A)= \frac{\area_g \{x\in Z, K(x)\in A\}}{\area_g Z}.
\end{equation}
We will write $\kappa(Z,A)$ when the metric is obvious.
Our first result asserts that statistically, a uniform part of the area of a large degree random curve is very curved. 
\begin{theorem}\label{syp} 
There exists $c>0$ such that  
	$$
\forall d\gg1,\  \mu_d\left[P\in \C^{hom}_d[X_0, X_1,X_2], \ 
\kappa\left(Z(P), g_{\mathrm{FS}|Z(P)}, [-4d,-\frac{d}8]\right) >  c\right]
\geq  c.$$
\end{theorem}
Recall that by the Gauss--Bonnet formula, the mean value of $K$ over a degree $d$ curve is $-d$. 
Note also that if we smooth $d$ different complex lines in $\C P^2$ passing through one common point, an arbitrary large proportion of the area of the resulting degree $d$ smooth curve have a curvature close to the one of a line, which implies that
\begin{equation}\nonumber
\forall \epsilon <2\pi, \ \forall d\geq 2,\   \inf_{P\in \C^{hom}_d[X_0, X_1,X_2]}
\kappa(Z(P), ]-\infty,\epsilon])
=0
\end{equation}
Hence, Theorem~\ref{syp} cannot be deduced from any deterministic fact. 
In a different direction, in~\cite{AG} it was proven that 
		$$  \mathbb \mu_d\left[
			\inf_{x\in Z(P)}              		K(x)\geq
		-d^{9}  \right]\underset{d\to \infty}{\to} 1.$$

Our second theorem computes exactly the average of the proportion of the area of the random complex curve with prescribed curvature. 
\begin{theorem}\label{deuxS}Let $0<r<R$. There exists $\varphi_{r,R}\in]0,1[$ defined below by~(\ref{mucha}) such that
	$$\forall d\geq 2, \ 
	\mathbb E_{\mu_d}\left[ 
	\kappa\Big(Z(P), [2\pi- Rd,2\pi-r d]\Big) 
	\right] =\varphi_{r,R},
	$$
	where $\kappa$ is defined above by~(\ref{kappa0}).
\end{theorem}
We now collect several consequences of Theorem~\ref{deuxS}.
\begin{corollary}\label{coro1}Let $0<r<R$. Then, for any $d\geq 2$,
	$$
	\begin{array}{rccl}
	\forall \eta \in ]\varphi_{r,R},1[, & 
	\mu_{d}\left[
	\kappa\Big(Z(P),  [2\pi-Rd,2\pi-r d]\Big) 
	>\eta 
	\right]& \leq &
	\frac{1}\eta \varphi_{r,R}\\
	\end{array}.$$
\end{corollary}
\begin{corollary}\label{coro2}Let $0<r<R$. Then, for any $d\geq 2$,
	$$
	\begin{array}{rccl}
	\forall \eta \in ]0,\varphi_{r,R}[, &
	\mu_d
	\left[
	\kappa\Big(Z(P),  [2\pi-Rd,2\pi-r d]\Big) 
	<\eta 
	\right]&\leq& 
	\frac{1}{1-\eta}(1- \varphi_{r,R}).
	\end{array}$$
\end{corollary}
Note that by the Gauss--Bonnet formula~(\ref{GB}),
for  any $r\geq 1$ and any $P\in \C^{hom}_d[X_0,X_1, X_2],$
$$
\kappa\Big(Z(P), ]-\infty,-rd]\Big) 
\leq \frac{1}{r}+O(\frac{1}d).
$$ 
\begin{corollary}\label{coro3}Let $\ell\in ]-\infty,2\pi[$. Then, 
	$$
	\begin{array}{rccl}
	\forall \eta\in ]0,1[,&  \mu_d \left[
	\kappa\Big(Z(P), [\ell,2\pi]\Big) 
>\eta
	\right]&\underset{d\to \infty}{\to} &0.
	\end{array}
	$$
\end{corollary}
The last assertion shows that, in particular, the event that  a uniform proportion of the area of the curve has a  non negative curvature becomes rarer and rarer for large degrees. We stress again that for any complex plane curve, there always exist points with positive curvature. 
 
\subsection{The general setting}\label{general}

 The previous results can be extended in a much more general setting that we now introduce.            
Let $S$ be a complex projective surface equipped with   a Hermitian ample holomorphic line bundle  $(L,h)\to S$  with positive curvature $\omega$, 
that is, locally 
$$ \omega = \frac{1}{2i\pi}\partial \bar{\partial}\log \|s\|^2_h>0,$$ 
where $s$ is any local non vanishing holomorphic section of $L$. Let $g_\omega= \omega (\cdot, i\cdot)$ be the associated K\"ahler metric.
The space $H^0(S, L^{d})$ of holomorphic sections of $L^{d}:=L^{\otimes d}$ is non trivial for $d$ large enough, more precisely $$ N_d:=\dim_\C H^0(S,L^d)\underset{d\to \infty}{\sim} d^2\vol_{g_\omega} (S). $$

Let $\Delta_d\subset H^0(S, L^d)$ be the discriminant subset, that is, the set of sections $s$ such that there exists $x$ in $Z(s)$ where $\nabla s(x)$ vanishes. Recall that $\Delta_d$ is a complex hypersurface, and that for any $s\in H^0(S, L^d)\setminus \Delta_d$, the zero set $Z(s)\subset S$ is a compact smooth complex curve of $S$. Moreover, since $H^0(S, L^d)\setminus \Delta_d$ is connected, for $s$ outside $\Delta_d$ the diffeomorphism class of $Z(s)$ depends only on $d$. For any $s$ we equip $Z(s)$ with the restriction $g_{\omega|Z(s)}$ of the K\"ahler metric $g_\omega$. 
By Wirtinger's theorem, $$\forall d\geq 1, \forall s\in H^0(S,L^d)\setminus \Delta_d, \ \area_{g_{\omega|Z(s)}} (Z(s))=\int_{Z(s)} \omega = 2d \vol_{g_\omega} (S).$$
The space $H^0(S,L^d)$ can be equipped with the $L^2$ Hermitian product 
\begin{equation}\label{prod}
(s,t)\in H^0(S,L^d)^2\mapsto \langle s,t\rangle = \int_S \langle s(x),t(x)\rangle_{h^d}\frac{\omega^2}{2}.
\end{equation}
This product induces a Gaussian measure $\mu_d$ over $H^0(S,L^d)$, 
that is, for any Borelian $U\subset H^0(S,L^d),$
\begin{equation}\label{mesure}
\mu_d (U)=\int_{s\in U} e^{-\frac12 \|s\|^2} \frac{\mathrm{d}s}{(2\pi)^{N_d}},
\end{equation}
where $\mathrm{d}s$ denotes the Lebesgue measure associated to the Hermitian product~(\ref{prod}).
If $(S_i)_{i\in \{1, \cdots, N_d\}}$
is an orthonormal basis of this space, then 
$$ s =\sum_{i=1}^{N_d} a_i S_i$$
follows the law $\mu_d$ if the random complexes $\sqrt 2 a_i$ are i.i.d standard complex Gaussians, 
that is, $\Re a_i$ and $\Im a_i$ are independent centered Gaussian variables with variance equal to $1/2$.
Note that for any event depending only on the vanishing locus $Z(s)$ of $s\in H^0(S, L^d)$, the probability measure $\mu_d$ can be replaced by the invariant measure over the unit sphere $\mathbb S H^0(S, L^d)$ for the product~(\ref{prod}), or equivalently the Fubini-Study measure on the linear system $\mathbb P H^0(S, L^d)$.

\begin{example}\label{bobmarley}
When $S= \C P^2$ and $(L, h)=(\mathcal O(1),h_{\mathrm{FS}})$ is the degree
	1 holomorphic line bundle equipped with the standard Fubini-Study metric, then the vector
	space $H^0(S, L^d)$ is isomorphic to the space 
	$ \C_d^{hom}[X_0 , X_1 , X_2 ]$ of  degree $d$ homogeneous
	polynomials in $3$ variables, and we recover the standard Fubini-Study measure given by~(\ref{mesure2}).
\end{example}

Then, $(Z(s), g_{\omega|Z(s)})$ can be seen as a fixed real surface with a random metric. The following theorem is the generalization of Theorem~\ref{syp} in this context.

\begin{theorem}\label{sy}
	Let $S$ be a compact smooth complex surface equipped with an ample holomorphic line bundle $(L, h) \to S$ endowed with a Hermitian metric $h$ with positive curvature $\omega$ and $g_\omega$ be the induced K\"ahler metric.
 Then, there exist a universal constant $c>0$ such that  
		$$
	\forall d\gg1,\  \mathbb \mu_d\left[
	s\in H^0(S,L^d), \ 
	\kappa\left(Z(s), g_{\omega|Z(s)}, [-4d,-\frac{d}8]\right) 
>  c
\right]
	\geq  c,$$
	where $\kappa$ is defined by~(\ref{kappa0}).
	\end{theorem}
The deterministic facts in the standard setting extend in this general one.
First,  
	$$ \exists C>0,\  \underset{d\to \infty}{\liminf}\inf_{s\in H^0(S,L^d)}
		\kappa\Big(Z(s), ]-\infty,- C]\Big)=0.$$
Indeed, let $d_0$ be such that
			$\dim_\C H^0(S,L^{d_0})\geq 2$ and $F\subset H^0(S,L^{d_0})\setminus \Delta_{d_0}$ be a compact subset. Then, there exists $C>0$ such that for any $s\in F$, the curvature of $Z(s)$ is larger than $-C$. Now for any integer $k\geq 1$, let $s_1, \cdots, s_k $ be $k$ sections in $F$ such that for $i\neq j$, $Z(s_i)$ intersects transversally with $Z(s_j)$. 
			Then, for any positive $\epsilon$,  we can find a small perturbation
			$s$ of $s_1\otimes \cdots \otimes s_k\in H^0(S,L^{kd_0})$ such that the area of $Z(s)$ where the curvature is less than $-2C$ is smaller than $\epsilon$. 

Second,
$$\forall d\gg 1, \ \inf_{s\in H^0(S,L^d)} \inf_{x\in Z(s)}K(x) = -\infty.$$ Indeed, let $Z$ be a degree $d$ curve with a nodal singularity at $x$. Then, for any $C>0$, there exists a smoothing of $C$ such that the curvature near $x$ is less than $-C$.

Third, by~\cite[Proposition 9.2]{kobayashi}, see also Theorem~\ref{coucou}, for any complex curve $Z$ in a complex K\"ahler manifold $(S,g)$, the curvature $K$ of $Z$  for $g_{|Z}$ is bounded above by the holomorphic sectional curvature of $S$, so that	 there exists $C>0$ such that 
$$K\leq C.$$  

Finally,
by Gauss--Bonnet, for any generic $s\in H^0(S,L^d),$ 
$$\int_{Z(s)} K (x)\mathrm{d}\vol(x) = - 2\vol(S) d^2 +O(d),$$ so that 
for any $r\geq 1$, 
$$
\kappa\Big(Z(s),]-\infty,-rd]\Big) 
\leq \frac{1}{r}+O(\frac{1}d).
$$ 

Also, observe that by~\cite{AG},
                	$$  \mathbb \mu_d\left[
                	s\in H^0(S,L^d),\ 
                	\inf_{x\in Z(s)}              		K(x)\geq
                	-d^{9}  \right]\to 1.$$

The following theorem generalizes Theorem~\ref{sy}:
                \begin{theorem}\label{deux}Under the hypotheses of Theorem~\ref{sy},
let $0<r<R$ and $\varphi_{r,R}\in]0,1[$ be defined below by~(\ref{mucha}). Then,  
                		$$
                		\mathbb E_{\mu_d}\left[
                		\kappa\Big(Z(s),[-Rd,-rd]\Big)
\right] 
\underset{d\to \infty}{\to} \varphi_{r,R},
                		$$
                		where $\kappa$ is defined by~(\ref{kappa0}).
                	\end{theorem}
                The following corollary generalizes Corollaries~\ref{coro1}, \ref{coro2} and \ref{coro3}:
\begin{corollary}\label{coro1S}
	Under the hypotheses of Theorem~\ref{sy},
let $0<r<R$. Then, 
$$
\begin{array}{rccl}
\forall \eta \in ]\varphi_{r,R},1[, & 
                \underset{d\to \infty}{\limsup}\,	\mu_d
                \left[
\kappa\Big(Z(s),[-Rd,-rd]\Big)
>\eta \right]& \leq &
                	\frac{1}\eta \varphi_{r,R}\\
                	\text{ and }
\forall \eta \in ]0,\varphi_{r,R}[, &
                                \underset{d\to \infty}{\limsup}\,	\mu_d
                                \left[
                	\kappa\Big(Z(s),[-Rd,-rd]\Big)
                	<\eta 
                	\right]&\leq& 
                	\frac{1}{1-\eta}(1- \varphi_{r,R}),
\end{array}$$
where $\varphi_{r, R}$ is defined by~(\ref{mucha}).
                	Moreover let $a,b\in \R$, $a<b$. Then,
                	$$
                	\begin{array}{rccl}
                	\forall \eta\in ]0,1[,&  \mu_d 
                	\left[
                	\kappa\Big(Z(s),[a,b]\Big)
                	>\eta
                	\right]&\underset{d\to \infty}{\to} &0.
                	\end{array}
                	$$
                \end{corollary}
We emphasize that our results in this general setting are universal, in the sense that the constants involved are universal and in particular do not depend on $S$. This is due to the fact that the Bergman kernel has itself a local universal behaviour, see appendix~\ref{bergman}.

\section{Proofs of the theorems}

\subsection{Deterministic preliminaries}

For any Riemannian surface $(Z,g)$ of finite area, any $A\subset \R$,
let 
\begin{equation}\label{kappa}
T(Z,g,A)= \{x\in Z, K_g(x)\in A\},
\end{equation}
where $K_g$ denotes the Gauss curvature.
Note that 
$$ \kappa(Z,g,A)= \frac{\area_g (T(Z,g,A))}{\area_g (Z)},$$
where $\kappa$ is defined by~(\ref{kappa0}).

\begin{lemma}\label{lemma1}There exists $f_0 :\C^2 \to \C$ a holomorphic function  vanishing transversally such that 
	$$
\area_{g_0}\left(T\Big(Z(f_0)\cap \mathbb B,g_0,[-2,-\frac14]\Big) \right)> 1,$$
where $g_0$ denotes the standard metric over $\C^2$ and $\mathbb B\subset \C^2$ is the standard unit ball.
\end{lemma}
\begin{preuve}
	Let $ f_0: \C^2 \to \C$ be the holomorphic function  defined $$\forall (z,w)\in \C^2, \ f_0(z,w)= zw-\frac14.$$ Then, $Z(f_0)\cap \frac{1}2\mathbb B \neq \emptyset$, so that
	$$ \area (Z(g)\cap \mathbb B )>0.$$ 
	By \cite[Proposition 1]{vitter}, for any  holomorphic function $f : \mathbb B \to \C$, 
	for any $x\in Z(f)$, the Gaussian curvature $K_{g_0}(x)$ of $Z(f)$ at $x$ for the standard metric equals 
	\begin{equation}\label{kzero}
	 K_{g_0}(x) =
	-\frac{|2f_{zw}f_zf_w-f_{zz}f_w^2 - f_{ww}f_z^2|^2}{(|f_z|^2+|f_w|^2 )^3}(x).
	\end{equation}
One can check that the curvature of $Z(f_0)\cap \mathbb B$ belongs to $[-2,-1/4]$, 
and that
	the area of $Z(f_0)\cap \mathbb B$ satisfies
	$$ \area_{g_0}(Z(f_0)\cap \mathbb B)= \int_{|z|^2+\frac1{16|z|^2}\leq 1}
	(1+\frac{1}{16|z|^4})\mathrm{d}z\wedge \mathrm{d}\bar{z}\geq 
	\int_{\frac{1}4\leq |z|^2\leq \frac{1}2}
\mathrm{d}z\wedge \mathrm{d}\bar{z}> 1.
$$
In conclusion, 
$ \area_{g_0}\{ x\in Z(f_0)\cap \mathbb B, -2<K_{g_0}(x)<-1/4\}> 1.$
\end{preuve}
In the sequel, we want to prove that the area of $T$  (defined by~(\ref{kappa})) is continuous in the $\mathscr C^2$-norm as a map of the ambient metric and the defining function of the curve $Z$.  For this, let us recall some general facts about the curvature. 
Let $Z$ be a submanifold of a Riemannian manifold $(S,g)$, $x\in Z$ and let
\begin{eqnarray}\label{sigma}
\sigma : T_xZ\times T_x Z &\to& N_x Z\\ \nonumber
(X,Y)& \mapsto & (\nabla_X Y)^\perp,
\end{eqnarray}
where $NZ\subset TS$ denotes the normal bundle over $Z$,  $\nabla$  the Levi-Civita connection associated to $g$ and $(\nabla_X Y)^\perp$ the $g$-orthogonal projection of $\nabla_X Y$ onto $N_xZ$.
\begin{theorem}\label{coucou}
		Let $Z$ be a submanifold of the Riemannian manifold $S$ and $x $ be in $Z$.  
	\begin{enumerate}
\item	\label{gauss}
	(Gauss's equations~\cite[Theorem 3.6.2]{jost2008riemannian})
	For any tangent vector $X,Y,V,W$ in $T_xZ,$
	\begin{eqnarray*}
		\langle R^Z(X,Y)V,W\rangle_g&=& \langle R^S(X,Y)V,W\rangle_g  +
		\\ &&\langle \sigma(Y,V),\sigma(X,W)\rangle_g- \langle \sigma(X,V),\sigma(Y,W)\rangle_g,
	\end{eqnarray*}
where $\sigma$ is defined by~(\ref{sigma}) and $R^Z$ (resp. $R^S$) denotes the Riemannian curvature of $g_{|Z}$ on $Z$ (resp. $g$ on $S$).
\item (K\"ahler version \cite[Proposition 9.2]{kobayashi}) Assume furthermore that $S$ is a K\"ahler manifold and  $Z$ is a complex submanifold of $S$. Then,
for any tangent vector $X \in T_xZ,$
\begin{eqnarray*}\label{kob}
\langle R^Z(X,JX)JX,X\rangle_g &=& \langle R^S(X,JX)JX,X\rangle_g  -\|\sigma(X,X)\|_g^2,
\end{eqnarray*}
where $J$ denotes the complex structure of $S$ and $g$ is the K\"ahler metric.
\end{enumerate}
\end{theorem}
In particular, as said in the introduction, 
the Gauss curvature of a complex curve is bounded by the holomorphic sectional curvature of its ambient space along the complex direction provided by the curve.

\begin{lemma}\label{lemma2} Let $0<r<R$ and $f_0 :2\mathbb B\subset  \R^4 \to \R^2$ be a $\mathscr   C^2$-function. Assume that $f_0$ vanishes transversely and that 
	$$\area_{g_0} \left(T\Big(Z(f_0),g_0,[-R,-r]\Big)\right)>0,$$
	where $k$ is defined by~(\ref{kappa}). Then, there exists $\delta>0$ such that for any $\mathscr   C^2$  metric $g$ and $\mathscr   C^2$ function $f : 2\mathbb B \to \R^2$ satisfying 
	$\|g-g_0\|_{\mathscr   C^2(2\mathbb B)}+\| f-f_0\|_{\mathscr   C^2(2\mathbb B)}\leq \delta,$ 
$$ \area_g\left(T(Z(f),g, [-2R,-\frac{r}2])\right)
>\frac{1}2\area_{g_0} \left(T\Big(Z(f_0),g_0,[-R,-r]\Big)\right).$$
\end{lemma}
\begin{preuve}By hypothesis, by the implicit function theorem and by compacity,  there exist $\eta>0$ and 
a finite set of cubes of the form $C_i=x_i+Q_\eta^1\times Q_\eta^2$, where $x_i\in \mathbb B$, $Q_\eta^i=]0,\eta[^2$, $i=1,2$, and $$ \mathbb B\subset \cup_i C_i \subset 2\mathbb B,$$ such that for any $\mathscr C^2$-function $f: 2\mathbb B \to \R^2$ $\mathscr C^1$-close enough to $f_0$,  for any $i$, $Z(f)\cap C_i $ is the graph of a $\mathscr   C^2$-function $\psi_i(f)$ over the translated of $Q_\eta^1$ or $Q_{\eta}^2$. Fix $i$ and assume that the graph $\psi_i(f)$ is defined over $Q_\eta^1$ in $C_i$. Let 
$$h_i(f):x\in Q_\eta^1\mapsto (x,\psi_i(f)(x))\in C_i.$$ 
By the implicit function theorem with parameters, $h_i(f)$ depends on $f$ continuously in the $\mathscr C^2$-norm. Then, for any $A\subset \R$, 
\begin{equation}\label{aire}
\area_g\{x\in Z(f)\cap C_i, K_g(x)\in A\}= 
\area_{h_i(f)^*g}
\{x\in Q_\eta^1, K_{g}(h_i(f)(x))\in A\}.
\end{equation}
In order to estimate the continuity of the right-hand side in $g$ and $f$, let us recall how to compute $K$ with respect to $g$ and $f$. 
Let $x\in Q_\eta^1$ and 
let 
\begin{eqnarray*}
	s: (T_x Q_\eta^1)^2& \to & N_{h_i(f)(x)}Z(f)\\
	(V,W)& \mapsto & \sigma (\mathrm{d}h_i(f)(x)(V),\mathrm{d}h_i(f)(x)(W)),
	\end{eqnarray*}
where $\sigma$ is defined by~(\ref{sigma}).
Recall that $\R^2\times \{0\}= T_x Q_\eta^1$
and $\mathrm{d}h_i(f)(x)(V)\in T_{h_i(f)(x)}Z(f).$
Then by~\cite[Proposition 3.2]{AG}, for any pair $(V,W)\in
(T_x Q_\eta^1)^2$,  
$$	s(V,W)= -
	\sum_{k=1}^2 \Phi^{-1}\left((\langle W_i,   \nabla^g_{V_i} \nabla f_j\rangle_g)_{j\in \{1, \cdots, 2\}}\rangle\right)_i \nabla f_k,
	$$
	where everything is computed at $h_i(f)(x)\in Z(f)$, where $\nabla f_k$ denotes the $g$-gradient of the $k$-th coordinate $f_k$  of $f$ and where 
			$$ \Phi = \left(\langle \nabla f_k, \nabla f_j\rangle_{g}\right)_{1\leq k,j\leq 2}.$$
	In particular, $s(V,W)$ is continuous in the $\mathscr   C^1$-norm in $g$, and in the $\mathscr   C^2$-norm in $f$.
Now by the Gauss equations (Theorem~\ref{coucou}), 
for any $x\in Q_\eta^1$ and 
and any $ X,Y,Z,W\in T_xQ_\eta^1 =\R^2\times \{0\},$
writing $X_i = \mathrm{d}h_i(f)(x)(X)\in T_{h_i(f)(x)}Z(f)$ etc., 
\begin{eqnarray*}
	\langle R^{Z(f)}(X_i,Y_i)Z_i,W_i\rangle_g&=& \langle R_g(X_i,Y_i)Z_i,W_i\rangle_g +
	\\ &&\langle s(Y,Z),s(X,W)\rangle_g- \langle s(X,Z),s(Y,W)\rangle_g,
\end{eqnarray*} 
where everything is computed at $h_i(f)(x).$
Now, the Riemannian curvature 
$ R_g$ depends
continuously of $g$ in the $\mathscr   C^2$-norm.
Consequently, if $K_{g_0}(h_i(f_0)(x))\in ]-R,-r[$, then uniformly
in $x\in Q_\eta^1$, for $g$ $\mathscr   C^2$-close to $g_0$ 
and $f$ $\mathscr   C^2$-close to $f_0$, then
$K(h_i(f)(x))\in ]-2R,-r/2[.$

Now, on $Q_\eta^1$ the pull-back metric $h_i(f)^*g$ converges  uniformly to $h_i(f_0)^*g_0$ when $g$ converges in $\mathscr   C^0$-norm to $g_0$ and $f$ to $f_0$ in the $\mathscr C^1$-norm,  so the quotient of their area form is bounded by 2 for $(g,f)$ close enough to $(g_0,f_0)$. Hence,
by~(\ref{aire}) the area of $T(Z(f)\cap C_i,g,r/2,2R)$ is larger or equal to half of the same for $g=g_0$ and $f=f_0$. Finally, using  a partition of the unity associated to the $C_i$'s, we obtain the result.
	\end{preuve}


\subsection{The Bargmann-Fock field}
The proofs of Theorems~\ref{sy} and~\ref{deux}
rely on similar probabilistic estimates 
for the universal algebraic rescaled model case, namely the Bargmann-Fock field over $\C^2$, see Theorem~\ref{Dai}. The Bargmann-Fock field is defined by the following measure:
\begin{equation}\label{bafo}
\forall z=(z_1,z_2)\in \C^2, \ f(z)=\sum_{(i_1,i_2)\in \Nn^2}a_{i_0, i_2}
{\sqrt{\frac{\pi^{i_1+i_2}}{i_1!i_2!}}}{z_1^{i_1} z_2^{i_2}}e^{-\frac12 \pi\|z\|^2},
\end{equation}
where the $a_I$'s are independent normal complex Gaussian random variables. We denote this measure by $\mu_{BF}$. note that, up to the exponential factor, its support is the set of holomorphic functions over $\C^2$. 

\subsection{Proof of Theorem~\ref{sy}}
Theorem~\ref{sy} is a  consequence
of the following local estimate:
\begin{proposition}\label{prop1}
	Under the hypotheses of Theorem~\ref{sy}, there exists a universal constant $c>0$ such that for any $d$ large enough and any $x\in S$, 
		$$ 
 \mu_d\left[
 s\in H^0(S,L^d),\ 
 	\area_{g_{\omega|Z(s)}}\left\lbrace y\in Z(s)\cap B(x,\frac{1}{\sqrt d}), -4d< K(y) < -\frac{d}8\right\rbrace> \frac{1}{2d} \right]\geq c.$$
\end{proposition}
We emphasize that the constant $c$ does not depend on $S$. 
\begin{preuve}
	Let $x\in S$ and $R>0$ be such that $2R$ is less than the radius of injectivity of $S$ at $x$. Note that since $S$ is compact, $R$ can be chosen independently of $x$. Then the exponential map based at $x$ induces a chart near $x$ with values in $B_{T_xS}(0,2R)$. We identify a point in $S$ with its coordinates. The tangent space $(T_xS,g_\omega)$ is identified with $(\C^2,g_0)$. 
	For any degree $d\geq 1$, let 
\begin{eqnarray*}
	 \psi_d : B(0,1)& \to  &B(0,1/\sqrt d),\\
	 y & \mapsto &y/\sqrt d.
	 \end{eqnarray*}
 For $d\geq 1$, let $$ g_d = d\psi_d^*g_{\omega}.$$
 Then $g_d $ converges to the standard metric $g_0$ in the $\mathscr   C^2$-topology. For any $d$, the function $\psi_d$ and the trivialization given by the parallel transport and explained in the appendix~\ref{bergman}
 provide a sequence of Gaussian complex
 functions $(f_d)_{d\geq 1}$ defined on $\mathbb B$ induced by the measure $\mu_d$ defined by~(\ref{mesure}).
  Let $f_0: \C^2\to \C$ be the map given by Lemma~\ref{lemma1}. Let $\delta>0$ given by Lemma~\ref{lemma2} applied to  $f_0$ and $(r,R)= (\frac{1}4,2)$ and
  $A$ be the event 
  $$ A=\left\lbrace f, \|f-f_0\|_{C^2(2\mathbb B)}< \delta/2\right\rbrace.$$
By Lemma~\ref{lemma2}, for $d$ large enough, 
  \begin{equation}\label{quoi}
   f\in A\Rightarrow
  \area_{g_d}\left(T\Big(Z(f),g_d,[-\frac{1}8,4]\Big)\right)>\frac12,
  \end{equation}
  where $T$ is defined by~(\ref{kappa}).
  Since by Theorem~\ref{Dai}, the kernel of $f_d$ converges in the $\mathscr   C^\infty$-topology to the kernel of the Bargmann-Fock field, 
by \cite[Theorem 4]{lerario2021differential}, 
$$  \liminf_{d\to \infty}
\mu_d [f_d\in A]\geq \mu_{BF}[f\in \inte A].$$
Now $\inte A = A$ and since $f_0$ lies in the support
of $\mu_{BF},$
there exists $c>0$ such that 
$$ \liminf_{d\to \infty}
\mu_d [f_d\in A]\geq  c.$$
By~(\ref{quoi}), this implies that
for $d$ large enough,
$$ \mu_d \left[\area\left(T\Big(Z(f),g_d,[-\frac{1}8,4]\Big)\right)>\frac12 \right]\geq  c,$$
hence the result after dilation by $1/\sqrt d$.
	\end{preuve}

\begin{preuve}[ of Theorem~\ref{sy}]
	For any $d\geq 1$, let $\Lambda_d$ be a maximal finite subset of points $x_i\in S$ such that the balls $B(x_i,\frac{1}{\sqrt d})$ are disjoint. Then, it is easy to see that for any $\epsilon >0$, 
	\begin{equation}\nonumber
	\forall d\gg 1,\  |\Lambda_d| \geq \frac1{16+\epsilon}\frac{\vol S}{\vol \mathbb B} d^2,
	\end{equation}
		see~\cite[2.5]{gayet2015expected}.
	Since $\vol \mathbb B =\pi^2/2,$
	\begin{equation}\label{yack}	
	\forall d\gg 1,\  |\Lambda_d| \geq \frac{\vol S}{80} d^2.
	\end{equation}
	For any $x\in \Lambda$, we define the following event $A(x)$:
	$$A(x)=\left\lbrace s\in H^0(S,L^d), \area\{y\in Z(s)\cap B(x,\frac{1}{\sqrt d}), -8d< K(y)<-\frac{d}8\} \geq \frac{1}{2d}\right\rbrace.$$
	By Proposition~\ref{prop1}, there exists a universal constant $c>0$ such that
	$$\forall d\gg1,  \ \forall x\in \Lambda_d, \ \mu_d [A(x)]\geq c.$$
	Let $N$ be the random variable defined by
	$$ N:=\# \{x\in \Lambda_d, A(x)\}.$$ Then, 
	\begin{eqnarray*}
		c|\Lambda_d|
		&\leq &
		\sum_{x\in \Lambda_d} \mu_d [A(x)] = \mathbb E \sum_{x\in \Lambda_d} {\bf 1}_{s\in A(x)}\\
		& \leq & \sum_{j=1}^{|\Lambda_d|}j \mu_d [N=j] \\
		&\leq &   \frac{c}2 |\Lambda_d| \mu_d
		\left[N\leq  \frac{c}2 |\Lambda_d|\right]
		+ |\Lambda_d |
		\mu_d \left[N\geq \frac{c}2|\Lambda_d| \right]
	\end{eqnarray*}
	which implies that
	$ 	\mu_d \left[N\geq \frac{c}2|\Lambda_d| \right]
	\geq \frac{c}2$
	and hence by~(\ref{yack}),
	$$ \mu_d\left[\area\{x\in Z(s), -4d<K(x)<-\frac{d}8\}
	)> c \frac{\vol S}{320}d \right]\geq \frac{c}2,$$
	hence the result.
\end{preuve}

\subsection{Proof of Theorem~\ref{deux}}

For Theorem~\ref{deux}, 
we will use the following Bargmann-Fock estimate:
\begin{proposition}\label{EBF} Let $U\subset \C^2$ be a bounded open subset with smooth boundary and $f: \C^2\to \C$ be the Bargmann-Fock field defined by~(\ref{bafo}). Then, for any $r,R>0$, 
	$$
	\mathbb E_{\mu_{BF}} \left[
	\area\{x\in Z(f)\cap U, -R<K(x)<-r\}
	\right] 
	= 2\vol (U)\varphi_{r,R}, 
	$$
	where 
	\begin{equation}\label{mucha} \varphi_{r,R} :=
	\int_{(a,b, \alpha, \beta, \gamma)\in \C^5}
	(|a|^2 +|b|^2)
	e^{-|(a,b,\alpha,\beta,\gamma)|^2/2}
	{\bf 1}_{\{\pi\frac{|2\gamma ab - \sqrt 2\alpha b^2 - \sqrt 2 \beta a^2|^2}{(|a|^2+|b|^2  )^3}\in [r,R]\}} \frac{\mathrm{d}(a,b,\alpha,\beta,\gamma)}{4(2\pi)^5}.
	\end{equation}
\end{proposition} 
Note that 
$$ \varphi_{0, \infty} :=
2\int_{a\in \C}|a|^2
e^{-|a|^2/2}
\frac{\mathrm{d}a}{8\pi}=1,$$ so that 
$\varphi_{r,R}\in ]0,1[.$
\begin{preuve}Let $0<r<R$. By the Kac--Rice formula, see for instance~\cite{azais2009level},
	$$ \mathbb E_{\mu_{BF}} [
	\area\{x\in Z(f)\cap U, -R<K(x)<-r\}
	]$$
	equals
	$$ 
	\int_U \mathbb E\left[{\bf 1}_{\{-R<K(x)<-r\}}|\mathrm{d}f(x)|^2 \, \big\vert \, f(x)=0
	\right]\rho_{f(x)}(0)\mathrm{d}x,$$
	where $\rho_{f(x)}(0)$ is
	the density of $f(x)$ at $0$, that is $(2\pi)^{-1}.$
	By invariance of the covariance function $\mathcal P$ of the BF field, see appendix~\ref{bergman}, 
	we can assume that $x=0$. Then, $f(0)$ and $\mathrm df(0)$ are independent, 
	and $$\forall (z,w)\in \C^2, \ f(z,w)= 
	\sqrt{\pi}(az+bw)+\pi\frac{\alpha}{\sqrt 2}z^2 +\pi\frac{\beta}{\sqrt 2}w^2 +\pi\gamma zw +O(3),$$
	where $a,b,\alpha, \beta $ and $\gamma$ are independent complex standard Gaussians,
	so that by~(\ref{kzero}) the curvature of $Z(f)$ at $(0,0)$
	equals $\pi\frac{|2 \gamma ab - \sqrt 2\alpha b^2 - \sqrt 2 \beta a^2|^2}{(|a|^2+|b|^2  )^3}$,
	hence the result. 
\end{preuve}

For the convenience's reader, we prove now Theorem~\ref{deuxS}, that is Theorem~\ref{deux} in the standard projective setting, which is easier and clearer than the general one. 
\begin{preuve}[ of Theorem~\ref{deuxS}]
Let $A\subset \R$ be a Borel subset. By the Kac--Rice formula, see for instance~\cite{azais2009level},
$$ \mathbb E 
\left[
\area\left\lbrace
x\in Z(P), K_{g_{\mathrm{FS}}}(x)\in A
\right\rbrace
\right] =
\int_{\C P^2} \mathbb E
\left[ {\bf 1}_{K_{g_{\mathrm{FS}}}(x)\in A}|\det{}^\perp \nabla P| \,\big\vert \, P(x)=0 
\right]\rho_{P(x)}(0 )
\mathrm{d}\vol (x),
$$
where $\rho_{P(x)}(0 )$ denotes the density at $0$ of the Gaussian field $P(x)$. 
By symmetry, we can assume that $x= [1:0\cdots : 0]\in \C P^2$. 
For $X_0\not=0$, let $f: \C^2\to \C$ be the random holomorphic function defined by
\begin{equation}\nonumber f:=\sqrt{\frac{2}{(d+2)!}}\frac{P}{X_0^d}
\end{equation}
Then, by~(\ref{mesure2}),
\begin{equation}\label{smoke}
 \forall (z,w)\in \C^2, \ f(z,w) = \sqrt d(az+bw)+ d(\frac{1}{\sqrt 2}\alpha z^2+\frac{1}{\sqrt 2}\beta w^2 +\gamma zw)+R_d,
 \end{equation}
where $R_d$ is a random polynomial vanishing at order 3 at $0$ and independent of the coefficients $a,b,\alpha, \beta$ and $\gamma$, which are standard complex Gaussians. 
If $d=1$, $Z(P)$ has constant curvature equal to $2\pi$. Assume that $d\geq 2$. 
By~\cite[p. 60]{ness1977curvature}, the curvature of $Z(f)$ at $0$ is
\begin{equation}\label{kz}
K_{g_{\mathrm{FS}}}(0)=
2\pi-\pi d\frac{|2\gamma ab - \sqrt 2\alpha b^2 - \sqrt 2 \beta a^2|^2}{(|a|^2+|b|^2  )^3}.
\end{equation}
The presence here of $\pi$ compared to \cite{ness1977curvature} is due to the difference choice of the Fubini-Study metric in her paper. 
By~(\ref{smoke}) and~(\ref{kz}), using the fact that at $[1:0:\cdots :0],$ $g_{\mathrm{FS}}= \sqrt{\pi} g_0 $, and 
assuming that $A=[2\pi-Rd,2-rd]$,
\begin{eqnarray}\label{Kac}\nonumber
\mathbb E
\left[ {\bf 1}_{K_{g_{\mathrm{FS}}}(x)\in A}|\det{}^\perp \nabla P| \,\big\vert \, P(x)=0 
\right]
&= &
d\int_{(a,b, \alpha, \beta, \gamma)\in \C^5}
\pi(|a|^2 +|b|^2)
e^{-|a,b,\alpha,\beta,\gamma|^2/2, }\times\\
\nonumber
&&\times {\bf 1}_{\{2\pi-\pi d\frac{|2\gamma ab - \sqrt 2\alpha b^2 - \sqrt 2 \beta a^2|^2}{(|a|^2+|b|^2  )^3}\in A\}} \frac{\mathrm{d}(a,b,\alpha,\beta,\gamma)}{(2\pi)^5}\\
&= &2\pi d\varphi_{r,R}.
\end{eqnarray}
Consequently, for any $r, R>0$, since $\vol(\C P^2) =\frac{1}2$, 
and $\rho_{f(0)}(0)=(2\pi)^{-1}$, one has
$$\forall d\geq 2, \ \frac{1}d\mathbb E
\left[\area\left\lbrace
x\in Z(P), K_{g_{\mathrm{FS}}}(x)\in [2\pi-Rd, 2\pi-rd]
\right\rbrace
\right] = \varphi_{r,R}\in ]0,1[,
$$
where $\varphi_{r,R}$ is defined by~(\ref{mucha}).
\end{preuve}
In the sequel, for any $r, R>0$,
let
$$ \area_{r,R}:=\area\left\lbrace
x\in Z(s), K_{g_{\mathrm{FS}}}(x)\in [2\pi-Rd, 2\pi-rd]
\right\rbrace.$$
\begin{preuve}[ of Corollary~\ref{coro1} and \ref{coro2}]
For any $\eta\in ]0,1[$, since $\area_{r, R}< d$, one has
$$
\eta  \mu_d(\area_{r,R}>\eta d)< \frac{1}d\mathbb E (\area_{r,R})
< \eta  \mu_d(\area_{r,R}<\eta d)+1-\mu_d(\area_{r,R}<\eta d),
$$
so that by Theorem~\ref{deuxS}, for any $\eta \in ]\varphi_{r,R},1[,$
$$ \mu_d[\area_{r,R}>\eta d]\leq 
\frac{1}{\eta} \varphi_{r,R}
$$
and for any $\eta \in ]0,\varphi_{r,R}[$,
$$ \mu_d[\area_{r,R}<\eta d]\leq 
\frac{1- \varphi_{r,R}}{1-\eta}.
$$
\end{preuve}
\begin{preuve}[ of Corollary~\ref{coro3}]
Let $\ell\in ]-\infty, 2\pi[$. Then, by~(\ref{kz})
$$K\in [\ell,2\pi]\Leftrightarrow
\frac{|2\gamma ab - \sqrt 2\alpha b^2 - \sqrt 2 \beta a^2|^2}{(|a|^2+|b|^2  )^3}\leq \frac{2\pi-\ell}{\pi d}.$$
When $d$ grows to infinity, the indicator function ${\bf 1}_{\{K\in [\ell,2\pi]\}}$  converges to the indicator of the hypersurface $$\{(a,b\alpha,\beta, \gamma)\in \C^5, \ 2\gamma ab - \sqrt 2\alpha b^2 - \sqrt 2 \beta a^2=0\}.$$ By the dominated convergence theorem, the Kac--Rice formula~(\ref{Kac}) implies that 
$$ \mathbb E_{\mu_d} [\area\{x\in Z(s),\ \ell< K(x)\leq 2\}]= o(d).$$
Now by Markov inequality, for any $\eta>0$,
$$ \mu_d[\area\{x\in Z(s),C< K(x)\leq 2\} > \eta d]= o(1),$$
hence the result.
\end{preuve}

In the general setting of random holomorphic sections of $H^0(S,L^d)$, we need to control more precisely the way the curvature of $Z(s)$ depends on $s$. For this, assume the hypotheses of Theorem~\ref{sy} and let $x\in S$. Let
$(t_i)_{i\in \{1, 2\}}$ be a real local orthonormal (for $h$) 
frame of $L$ given by parallel transport, 
and $G : TS^*\to TS$ be defined by 
\begin{eqnarray*}
	\forall \alpha\in TS^*,\  \langle G(\alpha),\cdot \rangle_g = \alpha,
\end{eqnarray*}
where $g:=g_\omega$ denotes the K\"ahler metric. Then~\cite[Proposition 3.2]{AG}, the bilinear operator $\sigma$ defined by~(\ref{sigma}) satisfies the formula:
for any $s\in H^0(S,L^d)\setminus \Delta_d,$
$$\forall V,W\in T_x Z(s), \  \sigma(V,W)=-
\sum_{i=1}^2 (\Phi^{-1}\left(\langle\nabla^2_{VW}s, t_j\rangle_h)_{j\in \{1, 2\}}\rangle\right)_i G \langle \nabla s,t_i\rangle_h,
$$
where 
	$ \Phi = \left(\langle \langle\nabla s,t_i\rangle_h, \langle\nabla s,t_j\rangle_h\rangle_{g^{*}}\right)_{1\leq i,j\leq 2}$,
and $g^*$ denotes the scalar product on $TS^*$ associated with $g$.
In the sequel, we will use the following homogeneity property of $\sigma$:
$$\forall \alpha\in \R, \beta\in \R^*, \  \sigma(V,W)( \nabla s, \nabla^2 s)
= \frac{\alpha}{\beta}\sigma(V,W)(\alpha \nabla s, \beta \nabla^2 s),$$
so that by~(\ref{kob}), 
\begin{eqnarray}\label{esmeralda}
	\ \forall X\in T_xZ(s),\ 
	\langle R^{Z(s)}(X,JX)JX,X\rangle( \nabla s, \nabla^2 s)&=& \langle R^S(X,JX)JX,X\rangle  -
	\\ \nonumber &&d\|\sigma(X,X)\|^2
( \frac{1}{\sqrt {d^{n+1}}}\nabla s, \frac{1}{\sqrt{d^{n+2}}}\nabla^2 s).
\end{eqnarray}

\begin{preuve}[ of Theorem~\ref{deux}]
Let $A\subset \R$ be a Borel subset. By the Kac--Rice formula, for any $d\geq 1$,
$$ \mathbb E 
\left[
\area\left\lbrace
x\in Z(s), K(x)\in A
\right\rbrace
\right] =
\int_{S} \mathbb E
\left[ {\bf 1}_{K(x)\in A}|\det{}^\perp \nabla s| \, \big\vert \, s(x)=0 
\right]\rho_{s(x)}(0 )
\mathrm{d}\vol (x),
$$
where $\rho_{s(x)}(0 )$ denotes the density at $0$ of the Gaussian field $s(x)$. 	
	Let $x\in S$. 
	Under the hypotheses and trivializations above near $x$ described in the appendix~\ref{bergman} 
	in any orthonormal basis of $T_xS$ (see for instance~\cite[Corollary 4.7]{bettigayet}):
	\begin{eqnarray}\label{coco} \cov\left(s, \nabla s, \nabla^2 s\right)_{|x} =
	d^2\left(	\begin{array}{ccc}
		(1+O(\frac1d))& O(\frac{1}{\sqrt d}) & O(1)\\
		O(\frac{1}{\sqrt d}) & \pi d{I}_{2}(1+O(\frac1d))&  O(\sqrt d)\\
		O(1)& O(\sqrt d) &  \pi^2 d^2\Sigma_{\mathrm{GOE}}(1+O(\frac1d))
	\end{array}\right),
	\end{eqnarray}
	where $I_2\in M_2(\R)$ is the identity matrix and $\Sigma_{\mathrm{GOE}}$ is defined by:
	\beq\label{goe}
	\Sigma_{\mathrm{GOE}} =  \left(\delta_{(ij)(kl)}+\delta_{(ji)(kl)}\right)_{\substack{1\leq i\leq j\leq 2\\
			1\leq k\leq l\leq 2}}\in M_{\frac{2(2+1)}2}(\C).
	\eeq
Let us define the random Gaussian variables 
$$ Q:= \frac{1}{\sqrt{d^2}}s(x),\ S := \frac{1}{\sqrt {d^{2+1}}}\nabla s(x) \text{ and } 
T:= \frac{1}{\sqrt{d^{2+2}}}\nabla^2 s(x).$$
Then, by~(\ref{coco}) and~\cite[Corollary 4.3]{bettigayet}, $\cov(R,S,T)$ converges to the covariance of
$(f(x), \mathrm{d}f(x), \mathrm{d}^{2}f(x))$ where $f$ is the Bargmann Fock field defined by~(\ref{bafo}). 

Now, by~(\ref{esmeralda}), for any $r, R>0$,
$$ K(s,\nabla s, \nabla^2 s)\in [-Rd,-r d]
\Leftrightarrow 
K(Q,S,T)\in [-R,-r] +O(\frac{1}d).
$$ 
Finally, when $d$ grows to infinity,
$$\frac{1}d \mathbb E
\left[ {\bf 1}_{K(x)\in [-Rd,r d]}|\det{}^\perp \nabla s(x)| \, \big\vert \, s(x)=0 
\right]\rho_{s(x)}(0 )$$ 
converges to
$$ \mathbb E_{\mu_{BF}}
\left[ {\bf 1}_{K_{g_0}(x)\in [-R,-r]}|\det{}^\perp \nabla f| \, \big\vert \, f(x)=0 
\right]\rho_{f(x)}(0 ).
$$
Hence, 
$$\frac{1}d \mathbb E 
\left[
\area\left\lbrace
x\in Z(s), K(x)\in A
\right\rbrace
\right] \to 2\vol (S) \varphi_{{r},{R}},
$$
where $\varphi$ is defined by~(\ref{mucha}).
By Wirtinger theorem, $\area(Z(s))= 2\vol(S)d,$
which concludes the proof of the theorem.
\end{preuve}
\begin{preuve}[ of Corollary~\ref{coro1S}]
The proof is the same as the one of the corollaries in the standard setting.
\end{preuve}

\appendix
\section{Asymptotics of the Bergman kernel}\label{bergman}

In this paragraph we assume that the setting and hypotheses of Theorem~\ref{sy} are satisfied. The covariance function $E_d$ for the Gaussian field generated by the holomorphic sections $s\in H^0(S, L^d)$ is defined by
$$ \forall z,w\in S, 
\ E_d (z,w) = \mathbb E \left[s(z)\otimes (s(w))^*\right]\in L^d_z \otimes (L^{d}_w)^*,$$
where the averaging is made for the measure $\mu_d$ given by~(\ref{mesure}), where $L^*$ is the (complex) dual of $L$ and 
$$\forall w\in S, \ \forall s,t\in L^d_w, \  s^* (t) = \langle s,t\rangle_{h^d(w)}.$$
The covariance $E_d$ is the \emph{Bergman kernel}, that is the kernel of the orthogonal projector from $L^2(M, L^d)$ onto $H^0(M, L^d)$. This fact can be seen through the equations
$$\forall z,w\in M, \ E_d(z,w) = \sum_{i=1}^{N_d} S_i (z)\otimes S_i^*(w),$$
where $(S_i)_i$ is an orthonormal basis of $H^0(M, L^d)$ for the Hermitian product~(\ref{prod}).
Recall that the metric $g_\omega$ is induced by the curvature form $\omega$  and the complex structure. 
It is now classical  that the Bergman kernel has a universal rescaled (at scale $\frac{1}{\sqrt d}$) limit, the Bargmann-Fock kernel $\mathcal P$:
\begin{equation}\label{mp}
\forall z,w\in \C^n, \ \mathcal P(z,w):=  \exp\left(
-\frac{\pi}2 (\|z\|^2 +\|w\|^2 -2\langle z,w\rangle)
\right).
\end{equation}
Theorem~\ref{Dai} below quantifies this phenomenon. For this, we need to introduce local trivializations and charts. 
Let $x\in S$ and $R>0$ be such that $2R$ is less than the radius of injectivity of $S$ at $x$. Then the exponential map based at $x$ induces a chart near $x$ with values in $B_{T_xS}(0,2R)$. We identify a point in $S$ with its coordinates. The parallel transport provides a trivialization 
$$\varphi_x : B_{T_xS}(0,2R)\times L^d_{x}\to  L^d_{|B_{T_xS}(0,2R) }$$ 
which induces a trivialization of 
$(L^d\boxtimes
( L^d)^*)_{|B_{T_xS}(0,2R)^2}$.
Under this trivialization, the Bergman kernel 
$ E_d$ becomes a map from $T_x M^2 $ with values into $\End\left( L^d_x\right)$.

\begin{theorem}(\cite[Theorem 1]{ma2013remark})\label{Dai} Under the hypotheses of Theorem~\ref{sy}, let $m\in \Nn$. Then, there exist $C>0$, such that for any $k\in\{0,\cdots, m\}, $ for any $x\in S$, 
	$\forall z,w\in B_{T_xS}(0,\frac{1}{\sqrt d}), $ 
	\beqe \left\|	\mathrm{d}^k_{(z,w)}\left(
	\frac{1}{d^n}E_d(z,w) - \mathcal P(z\sqrt d,w\sqrt d )
	\ \id_{L^d_x}
	\right)
	\right\|	
	\leq C d^{\frac{k}2-1}.
	\eeqe
\end{theorem}
The original reference is more general, see~\cite[Proposition 3.4]{letpuch} for the present simplification.\\

{\bf Acknowledgments.} The second author thanks Igor Wigman for a discussion about this subject. The research leading to these results has received funding from the French Agence nationale de la ANR-20-CE40-0017 (Adyct).

\bibliographystyle{amsplain}
\bibliography{courbure.bib}

\providecommand{\bysame}{\leavevmode\hbox to3em{\hrulefill}\thinspace}
\providecommand{\MR}{\relax\ifhmode\unskip\space\fi MR }
\providecommand{\MRhref}[2]{%
  \href{http://www.ams.org/mathscinet-getitem?mr=#1}{#2}
}
\providecommand{\href}[2]{#2}
\begin{thebibliography}{10}

\bibitem{AG}
Michele Ancona and Damien Gayet, \emph{Metric and spectral aspects of random
  complex divisors}, arXiv 2311.09679 (2023).

\bibitem{azais2009level}
Jean-Marc Aza{\"\i}s and Mario Wschebor, \emph{Level sets and extrema of random
  processes and fields}, John Wiley \& Sons, 2009.

\bibitem{bettigayet}
Damien Gayet, \emph{Expected local topology of random complex submanifolds}, J.
  Algebraic Geom. DOI: https://doi.org/10.1090/jag/817 (2022).

\bibitem{gayet2015expected}
Damien Gayet and Jean-Yves Welschinger, \emph{Expected topology of random real
  algebraic submanifolds}, Journal of the Institute of Mathematics of Jussieu
  \textbf{14} (2015), no.~4, 673--702.

\bibitem{jost2008riemannian}
J{\"u}rgen Jost, \emph{Riemannian geometry and geometric analysis}, vol. 42005,
  Springer, 2008.

\bibitem{kobayashi}
Shoshichi Kobayashi and Katsumi Nomizu, \emph{Foundations of differential
  geometry, volume 2}, Intersciences Publishers, 1969.

\bibitem{lerario2021differential}
Antonio Lerario and Michele Stecconi, \emph{Differential topology of {G}aussian
  random fields}, arXiv 1902.03805 (2021).

\bibitem{letpuch}
Thomas Letendre and Martin Puchol, \emph{Variance of the volume of random real
  algebraic submanifolds {II}}, Indiana Univ. Math. J. \textbf{68} (2019),
  no.~6, 1649–1720.

\bibitem{ma2013remark}
Xiaonan Ma and George Marinescu, \emph{{Remark on the off-diagonal expansion of
  the Bergman kernel on compact K{\"a}hler manifolds}}, Communications in
  Mathematics and Statistics \textbf{1} (2013), no.~1, 37--41.

\bibitem{ness1977curvature}
Linda Ness, \emph{Curvature on algebraic plane curves. {I}}, Compositio
  Mathematica \textbf{35} (1977), no.~1, 57--63.

\bibitem{SZ99}
Bernard Shiffman and Steve Zelditch, \emph{Distribution of zeros of random and
  quantum chaotic sections of positive line bundles}, Communications in
  Mathematical Physics \textbf{200} (1999), 661--683.

\bibitem{vitter}
Albert Vitter, \emph{On the curvature of complex hypersurfaces}, Indiana
  University Mathematics Journal \textbf{23} (1974), no.~9, 813--826.

\end{thebibliography}

\noindent
Michele Ancona \\
Laboratoire J.A. Dieudonn\'e\\
UMR CNRS 7351\\
Universit\'e C\^ote d'Azur, Parc Valrose\\
06108 Nice, Cedex 2, France\\

\noindent
Damien Gayet\\
 Univ. Grenoble Alpes, Institut Fourier \\
F-38000 Grenoble, France \\
CNRS UMR 5208  \\
CNRS, IF, F-38000 Grenoble, France
\end{document}